\documentclass{psp}

\newcounter{romanlistc}

\newenvironment{romanlist}
	{\setcounter{romanlistc}{0}
	 \begin{list}{$($\roman{romanlistc}$)$}
	{\usecounter{romanlistc}
	 \setlength{\parsep}{0pt}
	 \setlength{\itemsep}{3pt}
	 \setlength{\partopsep}{6pt}}}{\end{list}}

\usepackage{amssymb}
\usepackage{epsfig}

\newcommand{\Z}{{\mathbb Z}}
\newcommand{\Q}{{\mathbb Q}}
\newcommand{\Lk}{{\ell k}}
\newcommand{\VV}{{\mathcal V}}

\newtheorem{theorem}{Theorem}
\newtheorem{corollary}{Corollary}
\newtheorem{lemma}{Lemma}

\begin{document}

\title[The Conway potential function of a graph link]
  {The Conway potential function of a graph link}

\author[David Cimasoni]{DAVID CIMASONI\\
Section de Math\'ematiques, Universit\'e de Gen\`eve,\addressbreak
2-4, rue du Li\`evre, CP 240, CH-1211 Gen\`eve 24, Suisse\addressbreak
  e-mail\textup{: \texttt{David.Cimasoni@math.unige.ch}}}

\volume{136}
\pubyear{2004}
\setcounter{page}{1}
\receivedline{Received \textup{15} July \textup{2002}}
\maketitle

\begin{abstract}
We give a closed formula for the multivariable Conway potential function of any graph link in a homology sphere. As corollaries,
we answer three questions by Walter Neumann \cite{N} about graph links.
\end{abstract}

\section{Introduction}\label{sec:1}

A link $L$ in a homology $3$-sphere is called a {\it graph link\/} if its exterior is a graph manifold. For example, a link in $S^3$ is a graph link
if and only if it is {\it solvable}, that is, if it can be constructed by iterated cabling and connected sum operations from the unknot.
This class of links is particularly interesting, since all the links that arise in complex algebraic geometry are of this type.

Eisenbud and Neumann \cite{E-N} gave a classification of graph links by means of decorated trees called {\it splice diagrams\/}; they also found a formula
for the multivariable Alexander polynomial $\Delta_L$ of graph links in terms of the splice diagram (see (\ref{thm:multi-Alex}) below).
Since the Conway potential function is determined up to a sign by $\Delta_L$, the only problem is to compute this sign.

In 1999, Walter Neumann \cite{N} succeeded in computing the Conway polynomial (that is: the one variable potential function) of any fibered
solvable link. In this article, Neumann proves several formulas for fibered solvable links,
and asks whether these equalities still hold for any graph link. In this paper, we compute the multivariable potential function of any
(fibered or non-fibered) graph link in a homology sphere (Theorem \ref{theorem:CPF}). As a consequence, we give a positive answer to Neumann's
questions (Corollaries \ref{cor:1}, \ref{cor:2} and \ref{cor:3}). 

\section{The Conway Potential Function}\label{sec:2}

In 1970, Conway \cite{Con} introduced a new invariant of links called the {\it potential function}.
Given an oriented ordered link $L=L_1\cup\dots\cup L_n$ in $S^3$, its potential function is a well defined rational function
$\nabla_L(t_1,\dots,t_n)$ which satisfies
\begin{equation}\label{equ:Alex}
\nabla_L(t_1,\dots,t_n)\,\;\dot{=}\,\cases{\frac{1}{t_1-t_1^{-1}}\,\Delta_L(t_1^2) & if $n=1$; \cr           
            	\Delta_L(t_1^2,\dots,t_n^2) & if $n\geqslant 2$,\cr}
\end{equation}
where $\;\dot{=}\;$ means equal up to a multiplication by $\pm t_1^{\nu_1}\cdots t_n^{\nu_n}$. Thus, this invariant is basically the
multivariable Alexander polynomial without the ambiguity concerning multiplication by units of $\Z[t_1^{\pm 1},\dots,t_n^{\pm 1}]$. 
As a particular case of the potential function, Conway defined what was later called the {\it Conway polynomial} of a non-ordered link $L$.
It is given by
$$
\Omega_L(t)=(t-t^{-1})\;\nabla_L(t,\dots,t).
$$
Unfortunately, Conway's paper contains neither a precise definition of the potential function, nor a proof of its uniqueness.

In 1981, Kauffman \cite{Kau} found a very simple geometric construction of the Conway polynomial, namely
$$
\Omega_L(t)=\det\left(t^{-1}A-tA^T\right),
$$
where $A$ is any Seifert matrix of the link $L$ and $A^T$ the transpose of $A$. Finally, in 1983, Hartley \cite{Har} gave a definition of the
multivariable potential function $\nabla_L$ for any ordered oriented link in $S^3$. This definition was later extended by Turaev \cite{Tur} to
links in a $\Z$-homology $3$-sphere, and by Boyer and Lines \cite{B-L} to links in a $\Q$-homology $3$-sphere.

Let us now state several useful properties of the potential function; we refer to \cite{B-L} for the proofs.
Given an oriented ordered link $L=L_1\cup\dots\cup L_n$ in a $\Z$-homology sphere, there exists a well defined invariant $\nabla_L$
related to the multivariable Alexander polynomial $\Delta_L$ of $L$ by the equality (\ref{equ:Alex}). Furthermore, $\nabla_L$ satisfies the symmetry formula
\begin{equation}\label{equ:sym}
\nabla_L(t_1^{-1},\dots,t_n^{-1})=(-1)^n\,\nabla_L(t_1,\dots,t_n).
\end{equation}
Also, if $L'=(-L_1)\cup L_2\cup\dots\cup L_n$, where $(-L_1)$ denotes $L_1$ with the opposite orientation, then
\begin{equation}\label{equ:inv}
\nabla_{L'}(t_1,t_2,\dots,t_n)=-\nabla_L(t_1^{-1},t_2,\dots,t_n).
\end{equation}
Finally, if $L'=L_2\cup\dots\cup L_n$, we have the following {\it Torres formula\/}:
\begin{equation}\label{equ:Torres}
\nabla_L(1,t_2,\dots,t_n)=(t_2^{\ell_{12}}\cdots t_n^{\ell_{1n}}-t_2^{-\ell_{12}}\cdots t_n^{-\ell_{1n}})\;\nabla_{L'}(t_2,\dots,t_n),
\end{equation}
where $\ell_{ij}$ stands for the linking number $\Lk (L_i,L_j)$.

\section{Graph Links} \label{sec:3}
Eisenbud and Neumann \cite{E-N} gave a classification of graph links using splice diagrams. Following Neumann \cite{N}, we will
not recall the whole construction here, just several important features of these combinatorial objects. A splice diagram $\Gamma$ for a graph link $L$
is a tree decorated as follows:

\begin{romanlist}
\item{some of its valency one vertices are drawn as arrowheads, and correspond to the components of $L$;}
\item{the arrowheads have weight $+1$ or $-1$ (depending on whether the corresponding component of $L$ has ``intrinsic'' orientation or not);} 
\item{each edge has an integer weight at any end where it meets a {\it node} (vertex of valency greater than one), and all the edge weights
around a node are pairwise coprime.}
\end{romanlist}

Given a non-arrowhead vertex $v$ of $\Gamma$, there is a so-called {\it virtual component} $L_v\,$: this is the additional link component that would
correspond to a single arrow at $v$ with edge weight $1$. It is very easy, given two vertices $v$ and $w$, to compute the linking number $\Lk(L_v,L_w)$
of their corresponding components (virtual or genuine): if $\sigma_{vw}$ denotes the shortest path in $\Gamma$ joining $v$ and $w$ (including $v$
and $w$), then $\Lk(L_v,L_w)$ is the product of the edge weights adjacent to but not on $\sigma_{vw}$, and of the possible arrowhead weights.

Let us now recall Eisenbud and Neumann's formula for the multivariable Alexander polynomial of a graph link \cite[Theorem 12.1]{E-N}. 
If $L=L_1\cup\dots\cup L_n$ is a graph link given by a splice diagram $\Gamma$, its multivariable Alexander polynomial is equal to
\begin{equation}\label{thm:multi-Alex}
\Delta_L(t_1,\dots,t_n)\;\,\dot{=}\,\cases{(t_1-1)\prod_v(t_1^{\ell_{1v}}-1)^{\delta_v-2}& if $n=1$;\cr
				\prod_v(t_1^{\ell_{1v}}\cdots t_n^{\ell_{nv}}-1)^{\delta_v-2}& if $n\geqslant 2$,}
\end{equation}
where the product is over all non-arrowhead vertices $v$ of $\,\Gamma$, $\delta_v$ is the valency of the vertex $v$, and $\ell_{iv}$ denotes the
linking number of $L_i$ with $L_v$. In this equation, the terms of the form $(t_1^0\cdots t_n^0-1)^a$ should
be formally cancelled against each other before being set equal to zero.

In proving this theorem, Eisenbud and Neumann also show another remarkable result \cite[Theorem 12.2]{E-N}. Let us call an $n$-component link
{\it algebraically split\/} if, after possible renumbering, there is an index $1\leqslant q<n$ such that $\Lk(L_i,L_j)=0$ whenever $1\leqslant i\leqslant q<j\leqslant n$.
The theorem says that if $L$ is a graph link, then
\begin{equation}\label{thm:split}
\Delta_L(t_1,\dots,t_n)=0\quad\Longleftrightarrow\quad L\; \hbox{ is algebraically split}. 
\end{equation}

This is a very striking property of graph links. For example, it implies that the Alexander polynomial $\Delta_L$ of a $2$-component graph link is zero if
and only if the linking number $\ell$ of the components is zero. For general $2$-component links $L$, if $\Delta_L$ vanishes, then $\ell=0$ (by the Torres
formula). But the converse is false: the Whitehead link has Alexander polynomial $(t_1-1)(t_2-1)$, although $\ell=0$.
As a matter of fact, it is still an open question how to characterize geometrically $2$-component links with vanishing Alexander polynomial
(see \cite[Problem 16]{Fox}).

Before stating and proving our results, let us finally recall very briefly Neumann's argument \cite{N} for fibered solvable links.
Let $L$ be an oriented graph link given by a splice diagram $\Gamma$. Using the value of the one variable Alexander polynomial $\Delta_L$ of $L$, together
with the equality $\Omega_L(t)=\pm t^{-d}\Delta_L(t^2)$, where $d$ is the degree of $\Delta_L$, we obtain the formula
$$
\Omega_L(t)=\pm(t-t^{-1})\prod_v(t^{\ell_v}-t^{-\ell_v})^{\delta_v-2}\;,
$$
where the product is over all non-arrowhead vertices of $\Gamma$, $\delta_v$ is the valency of the vertex $v$, and $\ell_v$ denotes the linking
number of $L$ with $L_v$. Thus, the only issue is the determination of the sign of $\Omega_L$.

If $L$ is a fibered link, its Seifert matrix $A$ is unimodular. Therefore, the leading coefficient of $\Omega_L(t)=\det\left(t^{-1}A-tA^T\right)$
is given by $\det(-A)$. For fibered links, Lee Rudolph defined an integer invariant $\lambda$ called the
{\it enhanced Milnor number}. If the ambient sphere is $S^3$, this number is known to satisfy the following formula (see \cite{N-R}):
$$
(-1)^\lambda=\det(-A).\eqno(\star)
$$
Neumann proves that $\lambda\equiv k_-+j_-\!\!\!\pmod{2}$, where $k_-$ is the number of $(-1)$-weighted arrowheads and $j_-$
the number of non-arrowhead vertices $v$ for which $\ell_v$ is negative and $\delta_v$ is odd. This leads to
$$
\det(-A)=(-1)^{k_-+j_-},\eqno(\star\star)
$$
giving the formula
$$
\Omega_L(t)=(-1)^{k_-}(t-t^{-1})\prod_v(t^{\ell_v}-t^{-\ell_v})^{\delta_v-2}.\eqno(\star\star\star)
$$
These results lead Neumann to the following questions:

\begin{romanlist}
\item{is the formula $(\star\star\star)$ true for any graph link in a homology sphere~?}
\item{does the equality $(\star\star)$ hold for fibered graph links in a homology sphere~?}
\item{is the equation $(\star)$ still valid for fibered graph links in a homology sphere~?}
\end{romanlist}

Corollaries \ref{cor:1}, \ref{cor:2} and \ref{cor:3} provide an affirmative answer to these three questions.

\section{Results} \label{sec:4}

We are now ready to state and prove our results. Using (\ref{thm:multi-Alex}) along with the equations (\ref{equ:Alex}) and (\ref{equ:sym}),
it is easy to compute the Conway potential function up to a sign.

\begin{lemma}\label{lemma:uptosign}
The potential function of a graph link $L$ in a homology sphere is given by
$$
\nabla_L(t_1,\dots,t_n)=\epsilon(L)\;\prod_v(t_1^{\ell_{1v}}\cdots t_n^{\ell_{nv}}-t_1^{-\ell_{1v}}\cdots t_n^{-\ell_{nv}})^{\delta_v-2},
$$
where $\epsilon(L)$ is equal to $+1$ or $-1$.
\end{lemma}

\begin{proof*}
By (\ref{thm:multi-Alex}) and equation (\ref{equ:Alex}), $\nabla_L(t_1,\dots,t_n)$ is equal to
\begin{eqnarray*}
\nabla_L&=&\cases{\frac{\epsilon(L)}{t_1-t_1^{-1}}\; t_1^{\nu_1}\; (t_1-t_1^{-1})\;
		\prod_v(t_1^{\ell_{1v}}-t_1^{-\ell_{1v}})^{\delta_v-2} & if $n=1$; \cr           
            	\epsilon(L)\; t_1^{\nu_1}\cdots t_n^{\nu_n}\;
            	\prod_v(t_1^{\ell_{1v}}\cdots t_n^{\ell_{nv}}-t_1^{-\ell_{1v}}\cdots t_n^{-\ell_{nv}})^{\delta_v-2}  & if $n \geqslant 2$,}\\
&=&\epsilon(L)\; t_1^{\nu_1}\cdots t_n^{\nu_n}\;
		\prod_v(t_1^{\ell_{1v}}\cdots t_n^{\ell_{nv}}-t_1^{-\ell_{1v}}\cdots t_n^{-\ell_{nv}})^{\delta_v-2},
\end{eqnarray*}
for some integers $\nu_1,\dots,\nu_n$ and some sign $\epsilon(L)=\pm 1$.
The symmetry formula (\ref{equ:sym}) implies that $\nu_1=\dots=\nu_n=0$, giving the lemma.
\end{proof*}

\begin{lemma}\label{lemma:knot}
Let $K$ be a graph knot with $(+1)$-weighted arrowhead. Then $\epsilon(K)=+1$, that is:
$$
\nabla_K(t)=\prod_v(t^{\ell_v}-t^{-\ell_v})^{\delta_v-2}.
$$
\end{lemma}

\begin{proof*}
Let us use Kauffman's construction of the reduced potential function
$$
\Omega_K(t)=(t-t^{-1})\;\nabla_K(t)=\det(t^{-1}A-tA^T),
$$
where $A$ is any Seifert matrix for $K$. It is easy to check that the Seifert matrix for a graph knot $K$ is equal to the Seifert matrix for the
knot $K'$ obtained from $K$ by deleting every non-fibered splice component. (In other words: if $\Gamma$ is a splice diagram for $K$, a splice
diagram $\Gamma'$ for $K'$ is obtained from $\Gamma$ by deleting every vertex $v$ such that $\ell_v=0$.) Therefore, we can assume that $K$ is fibered.
In that case, $A$ is unimodular.
In particular, $\det(-A)$ is non-zero, so that $\det(-A)$ is the leading coefficient of $\Omega_L(t)=\det(t^{-1}A-tA^T)$.
There is an explicit computation of $A$ for graph knots in \cite{E-N}; it can be used to check the lemma.
Details can be found in \cite{these}.
\end{proof*}

\begin{lemma}\label{lemma:split}
Consider $L=L_1\cup\dots\cup L_n$. If $L-L_i$ is algebraically split for all $i$, then $L$ is algebraically split.
\end{lemma}

\begin{proof*}
Let us associate to $L$ a graph $G_L$ as follows: the vertices of $G_L$ correspond to the components of $L$, and two
vertices are linked with an edge if the linking number of the corresponding components is not equal to zero. Clearly, $L$ is algebraically split if and 
only if $G_L$ is not connected.
Given a vertex $v$ of a graph $G$, let us denote by $G-v$ the subgraph obtained by deleting the vertex $v$ and every edge adjacent to $v$.
We are left with the proof of the following assertion: given a graph $G$, if $G-v$ is not connected for any vertex $v$ of $G$, then $G$ is not connected.
In other words: if $G$ is a connected graph, there exists a vertex $v$ such that $G-v$ is connected. This last statement is very easy to prove:
given $G$ a connected graph, let $T$ be a maximal subtree of $G$. Since $T$ is a tree, it has at least
one vertex $v$ of degree one. Then, $T-v$ is connected, as well as $G-v$.
\end{proof*}

\begin{theorem}\label{theorem:CPF}
Let $L$ be a graph link with $n$ components given by a splice diagram $\Gamma$. Then, its Conway potential function is equal to
$$
\nabla_L(t_1,\dots,t_n)=(-1)^{k_-}\prod_v(t_1^{\ell_{1v}}\cdots t_n^{\ell_{nv}}-t_1^{-\ell_{1v}}\cdots t_n^{-\ell_{nv}})^{\delta_v-2},
$$
where the product is over all non-arrowhead vertices $v$ of $\,\Gamma$, $\delta_v$ is the valency of the vertex $v$, $\ell_{iv}=\Lk(L_i,L_v)$,
and $k_-$ is equal to the number of $(-1)$-weighted arrowheads.
\end{theorem}

\begin{proof*}
By formula (\ref{equ:inv}), it may be assumed that all the arrowheads have weight $(+1)$. Using the notation of Lemma \ref{lemma:uptosign},
we have to check that if $k_-=0$, then $\epsilon(L)=+1$. Let us do this by induction on $n\geqslant 1$.

The case $n=1$ is settled by Lemma \ref{lemma:knot}. Let us fix some $n\geqslant 2$, and assume that $\epsilon(L')=+1$ for all graph link $L'$ with
$n'\leqslant n-1$ arrowheads, all $(+1)$-weighted. Let $L=L_1\cup\dots\cup L_n$ be a graph link with $k_-=0$. If $\nabla_{L-L_i}=0$ for all $i$, then
(by (\ref{equ:Alex}) and (\ref{thm:split})) $L-L_i$ is algebraically split for all $i$. By Lemma \ref{lemma:split}, it follows that $L$ is
algebraically split, so that $\nabla_L=0$. In this case, there is nothing to prove. Therefore,
it may be assumed (after possible renumbering) that $\nabla_{L-L_1}\neq 0$. Let us note $L'=L-L_1$. If $\ell_{12}=\dots=\ell_{1n}=0$, then $L=L_1\sqcup L'$
is algebraically split, so $\nabla_L=0$. Hence, we can assume without loss of generality that $\ell_{12}\cdots\ell_{1n}\neq 0$. To summarize,
it may be assumed that
$$
(t_2^{\ell_{12}}\cdots t_n^{\ell_{1n}}-t_2^{-\ell_{12}}\cdots t_n^{-\ell_{1n}})\;\nabla_{L'}(t_2,\dots,t_n)\neq 0.
$$
The key ingredient of the induction step is the Torres formula (\ref{equ:Torres}). By Lemma \ref{lemma:uptosign}, we have
$$
\nabla_L(1,t_2,\dots,t_n)=\epsilon(L)\;\prod_{v\in \VV}(t_2^{\ell_{2v}}\cdots t_n^{\ell_{nv}}-t_2^{-\ell_{2v}}\cdots t_n^{-\ell_{nv}})^{\delta_v-2}.
$$
On the other hand, by (\ref{equ:Torres}),
$$
\nabla_L(1,t_2,\dots,t_n)=(t_2^{\ell_{12}}\cdots t_n^{\ell_{1n}}-t_2^{-\ell_{12}}\cdots t_n^{-\ell_{1n}})\;\nabla_{L'}(t_2,\dots,t_n)\;,
$$
which, by induction, is equal to 
$$
(t_2^{\ell_{12}}\cdots t_n^{\ell_{1n}}-t_2^{-\ell_{12}}\cdots t_n^{-\ell_{1n}})\;
\prod_{v\in \VV'}(t_2^{\ell_{2v}}\cdots t_n^{\ell_{nv}}-t_2^{-\ell_{2v}}\cdots t_n^{-\ell_{nv}})^{\delta_v-2}.
$$
Now, the set of non-arrowhead vertices $\VV'$ of $L'$ is equal to $\VV\cup \{v_1\}$, where $v_1$ is the vertex corresponding to the component $L_1$.
Therefore, $\delta_{v_1}=1$, $\ell_{v_1v}=\ell_{1v}$, and we have the equality
$$
\epsilon(L)\;\prod_{v\in \VV}(t_2^{\ell_{2v}}\cdots t_n^{\ell_{nv}}-t_2^{-\ell_{2v}}\cdots t_n^{-\ell_{nv}})^{\delta_v-2}=
\prod_{v\in \VV}(t_2^{\ell_{2v}}\cdots t_n^{\ell_{nv}}-t_2^{-\ell_{2v}}\cdots t_n^{-\ell_{nv}})^{\delta_v-2}.
$$
Since
$$
\prod_{v\in \VV}(t_2^{\ell_{2v}}\cdots t_n^{\ell_{nv}}-t_2^{-\ell_{2v}}\cdots t_n^{-\ell_{nv}})^{\delta_v-2}=
	(t_2^{\ell_{12}}\cdots t_n^{\ell_{1n}}-t_2^{-\ell_{12}}\cdots t_n^{-\ell_{1n}})\;\nabla_{L'}(t_2,\dots,t_n),
$$
which is assumed to be non-zero, $\epsilon(L)$ is equal to $+1$.
\end{proof*}

\begin{corollary}\label{cor:1}
The Conway polynomial of a fibered graph link $L$ is given by
$$
\Omega_L(t)=(-1)^{k_-}(t-t^{-1})\prod_v(t^{\ell_v}-t^{-\ell_v})^{\delta_v-2}.
$$
\end{corollary}

\begin{proof*}
Use the fact that $\Omega_L(t)=(t-t^{-1})\;\nabla_L(t,\dots,t)$.
\end{proof*}

\begin{corollary}\label{cor:2}
Let $A$ be a Seifert matrix for a fibered graph link $L$; then 
$$
\det(-A)=(-1)^{k_-+j_-},
$$
where $j_-$ is the number of non-arrowhead vertices $v$ with $\ell_v$ negative and $\delta_v$ odd.
\end{corollary}

\begin{proof*}
If $L$ is a fibered link, $A$ is unimodular, so $\det(-A)$ is the leading coefficient of
$\Omega_L(t)=\det(t^{-1}A-tA^T)$. By Corollary \ref{cor:1}, $\det(-A)$ is the leading coefficient of
$(-1)^{k_-}\prod_v(t^{\ell_v}-t^{-\ell_v})^{\delta_v-2}$, which is clearly $(-1)^{k_-+j_-}$.
\end{proof*}

\begin{corollary}\label{cor:3}
Let $A$ be a Seifert matrix for a fibered graph link $L$; then
$$
\det(-A)=(-1)^\lambda\,,
$$
where $\lambda$ is the enhanced Milnor number of $L$.
\end{corollary}

\begin{proof*}
By Theorem 6.1 and $\S$10 in \cite{N-R}, we have the equality
$$
\lambda=\ell k(P\cup L_+,N\cup L_-),
$$  
where the notations are as in \cite{N}. By Corollary \ref{cor:2}, it remains to show that
$$
\ell k(P\cup L_+,N\cup L_-)\equiv k_-+j_-\!\!\!\pmod{2}.
$$
The argument is exactly as in \cite{N}.
\end{proof*}

Of course, Corollaries \ref{cor:2} and \ref{cor:3} are false for non-fibered graph links. For example, a Seifert matrix $A$
for the trivial $2$-component link satisfies $\det(-A)=0$.

\begin{ex*}
{\rm Consider the graph link $L$ given by the splice diagram illustrated below,
 
\begin{figure}[ht] 
\begin{center}
\epsfig{figure=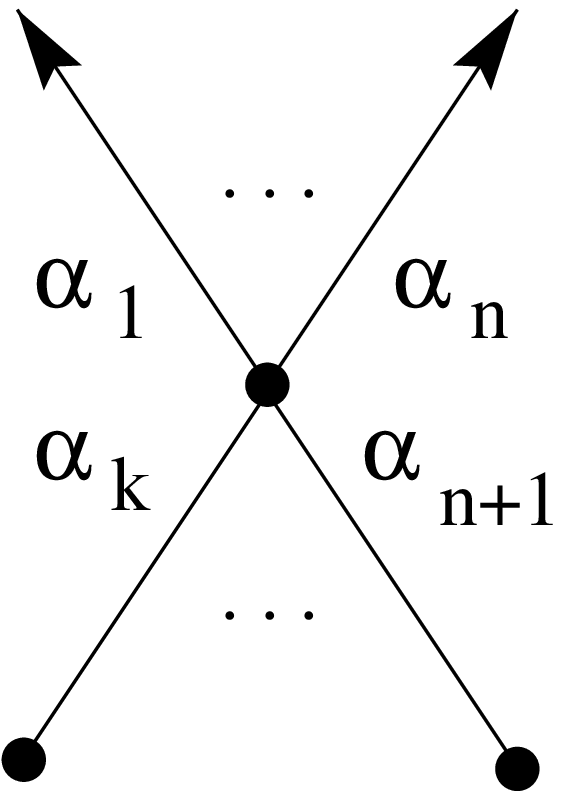,height=2.7cm}
\end{center}
\end{figure}
and let us note $\alpha=\alpha_1\cdots\alpha_k$. Then, the potential function of $L$ is equal to}
$$
\nabla_L(t_1,\dots,t_n)\;=\;\frac{\left(t_1^{\alpha/\alpha_1}\cdots t_n^{\alpha/\alpha_n}-t_1^{-\alpha/\alpha_1}\cdots t_n^{-\alpha/\alpha_n}\right)^{k-2}}
{\prod_{i=n+1}^k\left(t_1^{\alpha/\alpha_1\alpha_i}\cdots t_n^{\alpha/\alpha_n\alpha_i}-t_1^{-\alpha/\alpha_1\alpha_i}\cdots t_n^{-\alpha/\alpha_n\alpha_i}\right)}.
$$
\end{ex*}

\begin{acknowledgments}
The author wishes to express his thanks to Fran\c coise Michel, Claude Weber and Mathieu Baillif.
\end{acknowledgments}

\end{document}